\documentclass[12pt]{article}
\usepackage{geometry}                
\usepackage{graphicx}
\usepackage{amssymb}					
\usepackage{amsmath}					
\usepackage{mathrsfs}
\usepackage{tikz}
\usepackage{epstopdf}
\usepackage{amsthm}
\newtheorem*{theorem*}{Theorem}
\newtheorem*{lemma*}{Lemma}
\newtheorem*{example*}{Example}
\newtheorem{theorem}{Theorem}
\newtheorem{lemma}{Lemma}
\newtheorem{corollary}{Corollary}

\def\eop{\hfill $\blacksquare$ \smallskip}
\def\bop{\noindent {\bf Proof.} }

\def\b0{{\bf 0}}
\def\b1{{\bf 1}}

\def\n{\noindent}

\begin{document}
\title{Book embeddings of graphs and a
theorem of Whitney\thanks{see also https://faculty.georgetown.edu/kainen/pbip3.pdf
}}

\author{
Paul C. Kainen\\
 \texttt{kainen@georgetown.edu}
\\
Shannon Overbay
\\
\texttt{overbay@gonzaga.edu}
}
\date{}                                           

\newcommand{\Addresses}{{
  \bigskip
  \footnotesize

\n
Paul C. Kainen, \textsc{Department of Mathematics and Statistics,\\ Georgetown University,
37th and O Streets, N.W., Washington DC 20057}\\
\\

\vspace{-0.27cm}
\n

\par\nopagebreak
}}

\maketitle

\abstract{
It is shown that the number of pages required for a book embedding of a graph is the maximum of the numbers needed for any of the
maximal nonseparable subgraphs and that a plane graph in which every triangle bounds a face has a two-page book embedding. The latter
extends a theorem of H. Whitney and gives two-page book embeddings
for $X$-trees and square grids.
}

\smallskip

\n
{\bf Key Phrases}: {\it Hamiltonian planar graph, book thickness, nicely planar
graph, girth, bipartite.} This is the full version of  Extension of a theorem of Whitney, {\it Appl Math Lett.} 20 (2007) 835--837 by the same authors.

\vspace{0.5 cm}
\n

\section{Introduction}

In this paper, we consider finite simple graphs (no loops or parallel edges);
see, e.g.,\cite{harary}, \cite{dbw}. Let $bdy(e) = \{v, w\}$ denote the boundary vertices (endpoints) of an edge $e$.

By a closed half-plane, we mean a copy of the complex numbers with
imaginary part $\geq 0$.   For $k$ a positive integer, the {\bf book} $B_k$ is the singular
surface formed from the union of $k$ closed half-planes (the pages) intersecting in a line L (the spine) which is the (topological) boundary of each of the pages.
 
In \cite{PCKa} (p. 97), a {\bf k-page book embedding} of a graph $G = (V, E)$ is defined to be an
embedding of $G$ into $B_k$ which carries $V$ to $L$ with the property that each
edge $e \in E$ is mapped into a single page so that $e \cap L =bdy(e)$. 
The {\bf book
thickness}, or {\bf pagenumber}, of $G$, $bt(G)$, is the least number of pages in which
$G$ has a book embedding.

Equivalently, $bt(G)$ is the least number of colors
which suffice to color the edges of $G$ so that no two edges of the same color
intersect other than at a common endpoint when the vertices of $G$ are arranged
in order around the boundary of a circle and the edges are drawn as straightline segments, minimized over all possible cyclic orderings of the vertices [3].
It is easy to see that $G$ has book thickness 1 if and only if $G$ is outerplanar
and $bt(G) \leq 2$ if and only if $G$ is a subgraph of a planar Hamiltonian graph
\cite{bk79}.

Book thickness has turned up in rather diverse applications such as fault-tolerant computing and VLSI (Chung, Leighton, and Rosenberg \cite{clr1987}), computational complexity and graph separators (Galil, Kannan, and Szemeredi
\cite{gal1989}), software complexity metrics and vehicle traffic engineering (Kainen \cite{PCKb}),
and ``bisecondary structures'' used for modeling RNA folding energy states
(Gleiss and Stadler \cite{gle1999}). In addition, many theoretical questions remain
open such as the relationship of book thickness to other invariants \cite{dhutch1991}, \cite{PCKb},
and the book thickness of standard and useful families of graphs (e.g., \cite{MWW},
\cite{Heath}, \cite{Has1997}, \cite{eno1999}).

Recall that a graph is planar if it can be embedded in the plane; the
actual embedding is called a plane graph. We call a graph {\bf nicely planar} if it
has a plane embedding in which each triangle bounds a region.

In this paper, we investigate some connections between book thickness
and planarity. Yannakakis has shown \cite{Yann} that the book thickness of a planar
graph cannot exceed 4 and examples are known of maximal planar graphs
with no Hamiltonian cycle so 3 (and perhaps 4) pages are certainly necessary
for some planar graphs.  (As of 2021, it is known that 4 pages are necessary.)

Our main result is that a nicely planar graph is a subgraph of a Hamiltonian planar graph. We then apply the result to give simple proofs that
some useful graphs in computer science have two-page book embeddings,
previously shown by direct computations \cite{clr1987}. The proof of our main result
requires a lemma of independent interest: $bt(G) = max bt(B)$, where $B$ is a
maximal nonseparable subgraph of G.

In section 2, we study book thickness in terms of the block-cutpoint tree.  Section 3 gives the main result, while section 4 has applications.
 \section{ Book thickness and block decomposition}
 
A graph is {\it nonseparable} if it has no cutpoints; a {\it block} is a maximal nonseparable subgraph \cite{harary} (p. 26). A family of subsets of some given set determines
an {\it intersection graph} where vertices correspond to subsets and adjacency
holds if and only if the corresponding subsets (are distinct and) have nonempty intersection.
The intersection graph formed by the cutpoint singletons and the vertex sets
of the blocks is a forest, called the {\it block-cutpoint forest} \cite{harary} (p. 36). The
points of the intersection graph which have degree 1 must be blocks since
every cutpoint belongs to more than one block. A graph is connected if and
only if its block-cutpoint forest is a tree.

\begin{theorem}
The book thickness of a graph is the maximum of the book
thicknesses of its blocks.
\end{theorem}
\bop
Let G be a graph. Without loss of generality we assume G is connected. We prove the result by induction on the number of blocks. The basis
case of one block is trivial.

Consider the block-cutpoint tree of $G$ which by assumption is not trivial.
There is a degree-one vertex corresponding to a block $H$ with a unique
cutpoint $v$ and $G' = G \setminus (H - v)$ has one fewer block than $G$. By the
inductive hypothesis, $G'$ has book thickness equal to the maximum of that
of its blocks. So we now need only show that $G$ requires no more than
the maximum of the book thicknesses of $H$ and $G'$. Take minimum book
thickness embeddings of $H$ and $G'$ and, by rotating the vertices along the
spine if necessary, make $v$ the first vertex of each embedding and place the
two embeddings consecutively along the spine of a book, using the maximum
number of pages in the two embeddings - WLOG, take $G'$
second. Now pull
the copy of $v$ in the second embedding above $H$ and superimpose it on the
first copy of $v$. Hence, $bt(G) = \max\{bt(G'), bt(H)\}$, as required. \eop

Note that the genus of a graph is the sum of the genera of its blocks
(Battle, Harary, Kodama, and Youngs \cite{bhky1982}). For book thickness, maximum
replaces sum. With respect to block decomposition, {\it genus is to book thickness as the $\ell_1$-norm is to the $\ell_\infty$-norm}.

This type of distinction separates invariants, such as the cyclomatic number, which are measures of global complexity (additive over the blocks) from
invariants, such as clique size, which are local measures of complexity and
have values obtained by maximizing over the blocks.

We wonder whether there are interesting graph theoretic complexity measures corresponding to the $\ell_2$ (Hilbert) norm. Note that the global complexity
measures given here do not distinguish homeomorphic graphs, while those of
local type may have very different values for homeomorphic graphs; see the
last section

\section{Extending Whitney’s theorem}

A {\it triangulation} is a planar graph with a maximal set of edges. Equivalently,
a graph is a triangulation if it is isomorphic to a plane graph in which every
face, including the face which contains infinity, has exactly three edges in its
boundary. A {\it triangle} in a graph is a cycle of length 3.

By the Jordan Curve Theorem, in a plane graph, any triangle divides the
plane into an interior and an exterior region. A triangle is bounding if either
its interior or exterior region contains no vertices from the graph. A triangle
is separating if its deletion increases the number of connected components of
the graph. In a triangulation, a triangle is separating if and only if it is not
bounding. If there are no separating triangles in a planar graph, then the
graph is nicely planar, and the converse holds if the graph is 2-connected.

As an example, consider the plane graph $G$ which consists of two connected components: two disjoint triangles side by side with three additional
edges all attached like a bonnet at the top vertex of the first triangle, pointing up into the outer region of the first triangle. Removing the first triangle
increases the number of connected components from 2 to 4, the components
now being the second triangle and three isolated vertices. The first triangle
is separating but also bounding and the second triangle also bounds so $G$ is
nicely plane.

A graph is {\it Hamiltonian} if it has a cycle through all of the vertices. The following theorem is due to H. Whitney \cite{Whit}.
\begin{theorem}
Every triangulation with no separating triangles is Hamiltonian.
\end{theorem}

A graph is {\it subhamiltonian} if it is a subgraph of a planar Hamiltonian graph.
A plane graph $G'$
is said to be obtained by {\it stellating} a face $F$ of the plane
graph $G$ if $G'$
is obtained from $G$ by adding a new vertex $F^*$
in the interior
of $F$ and continuous curves joining it to each vertex in the boundary of the
face in such a way that the curves intersect $G$ only at their endpoints and
the curves are disjoint except for $F^*$.
See, e.g., \cite{bk79}.

\begin{theorem}
Every 3-connected planar graph with no separating triangles
is subhamiltonian.
\end{theorem}
\bop
Choose a plane embedding of the given graph $G$ and extend it by
stellating any nontriangular faces. The resulting graph $G'$
is clearly a planar triangulation and we claim that it has no separating triangles. Indeed,
suppose $T$ were a separating triangle in $G'$. Then at least one of the vertices
is new since $G$ had no separating triangles. But the only new vertices were
added as stellation points so $T = \{v, x, y\}$, where $v$ is a stellation point and
$x,y$ are vertices in the boundary of the face which was stellated by $v$. It is
straightforward to verify that $\{x, y\}$ would be a separating set for $G$ so no
separating triangles exist in $G'$
\cite{Overbay}. By Theorem 2, $G'$ and so
$G$ are subhamiltonian. 
\eop

\begin{lemma}
Let $G$ be a nicely planar block. Then there is a 3-connected
planar graph $G'$ with no separating triangles such that $G$ is a subgraph of $G'$.
\end{lemma}
\bop We induct on the number of separating 2-sets among the vertices of
$G$. If there are none, then $G$ is already 3-connected so the basis case of the
induction holds with $G' = G$.

Suppose $G$ has $k \geq 1$ separating 2-sets and take $A = \{u, v\}$ to be one
of these separating 2-sets. Choose a particular plane embedding of $G$. Let
$G - A = G_1 \cup \cdots \cup G_n$ with $G_j$ denoting the connected components obtained
by removing $A$, listed in the clockwise order of the edges joining them to $v$.
If $vu$ is an edge, we may renumber the components so that the edge from v
to u precedes all the edges from $v$ to $G_1$, and after all those from $v$ to $G_n$.
With respect to the same fixed clockwise order at $v$ there is a last vertex $v_j$
in $G_j$ and a first vertex $w_j$
in $G_{j+1}$, for $j = 1, \ldots, n-1$. Add $n-1$ vertices $y_j$
together with the edges $y_jv_j$, $y_jw_j$, and $y_jv$, keeping the resulting graph $G'$
plane. Indeed, $G'$
is nicely plane since the added triangles all bound disks.

Hence, it suffices to show that (1) $G'$
is 2-connected and (2) $G'$ has fewer
separating 2-sets than $G$. The assertion (1) follows from the ``ear'' decomposition characterization, e.g., \cite{dbw} (p. 163). One can form $G'$ by adding paths
(ears) in the order $v_1, y_1, w_1, vy_1, v_2, y_2, w_2, \ldots, vy_n$.

Assertion (2) follows from the claims: (i) No added vertex is a member
of a separating 2-set for $G'$,  (ii) $A = \{u, v\}$ is not a separating 2-set for $G'$, and (iii) if $B$ is any separating 2-set for $G'$, then $B$ is a separating 2-set for $G$.

These may be established as follows: For (i), first note that for $B$ any
separating 2-set of $G'$, exactly one of the vertices must be new. Since the old
vertex can’t be a cutpoint of the block $G$, one component of $G' - B$ consists of
a single new vertex. But this is impossible since each new vertex is adjacent
to 3 old ones. For (ii), let $A = \{u, v\}$ as given above in the definition of $G'$.
By our construction, $G' - A$ is connected so $A$ is not a separating set for $G'$.  Last, for (iii), let $B$ be a separating 2-set for $G'$. No component of $G' - B$
is a single new vertex by the same degree argument as above. Hence, each
component contains at least one vertex of $G$, so $B$ is a separating 2-set for
$G$. See \cite{Overbay} for a similar and slightly more detailed argument. \eop

\begin{theorem}
Every nicely plane graph is subhamiltonian.
\end{theorem}
\bop By Lemma 1 and Theorem 3, each block is subhamiltonian and
this suffices by Theorem 1. \eop

Note that we cannot always add edges alone to a nicely plane graph in a
way to satisfy the conditions of Whitney’s theorem. For example, a square
with one stellated face is nicely plane but extending to a plane graph by
adding one more edge forces a separating triangle.

Our method adds both vertices and edges to avoid creation of separating
triangles - hence, enabling the use Whitney’s Theorem to obtain a Hamiltonian cycle. However, once such a cycle is created, one can delete both the
added edges and vertices of the construction, keeping the ordering for the
original vertices. By adding edges joining any nonadjacent pairs of consecutive vertices, a Hamiltonian planar graph is obtained which contains the
original nicely planar graph as a spanning subgraph. The same argument
shows that in the definition of subhamiltonian graph, one can require that
the extension only involve the inclusion of new edges.

\section{Applications} 

An {\bf X-tree} is the plane graph formed by taking a complete plane binary tree
oriented with the root at the top (with branching down and to left and right,
resp.) and adding to it horizontal paths which join all vertices at the same
distance to the root. Let us call an {\bf extended X-tree} the result of using an
additional edge for each of the added horizontal paths, joining the endpoints
of the path to form a cycle unless these two vertices were already adjacent.
The following is an immediate consequence of Theorem 4 since an extended
$X$-tree is nicely plane; cf. Chung, Leighton and Rosenberg \cite{clr1987}
\begin{corollary}
Every extended X-tree (hence every X-tree) is subhamiltonian.
\end{corollary}
In \cite{clr1987} there is also an explicit proof that the product of two paths (called
a {\it square grid}) is subhamiltonian. This follows from the next corollary.

Recall that the {\it girth} of a graph is the length of its shortest cycle. Graphs
with girth $> 3$ are triangle-free and so if they are planar, then they are
nicely planar. In particular, {\it every bipartite planar graph has a 2-page book
embedding}.

\begin{corollary}
Every planar graph with girth $> 3$ is subhamiltonian.
\end{corollary}

Theorem 4 has implications for two open questions. Barnette asked if
every cubic 3-connected planar bipartite graph is Hamiltonian; see, e.g., \cite{Hol}.
By Corollary 2 such a graph is at least subhamiltonian. Chartrand, Geller
and Hedetniemi \cite{char1971} conjectured that every planar graph can be written as
the edge-disjoint union of two outerplanar graphs. By Theorem 4, their
conjecture is true for nicely planar graphs (Goncalves, STOC'05, claims a general proof).

Recall that two graphs are {\it homeomorphic} if they have isomorphic subdivisions.
It is well known that every graph is homeomorphic to a graph of book
thickness at most three (see Atneosen \cite{gat1968}, Bernhart and Kainen \cite{bk79}, and,
according to Jozef Przytycki, also a dissertation by G. Hotz, a student of
Reidemeister). Applications of three-page embeddings have been made to
links by Dynnikov \cite{dynn1999}.

In contrast, for planar graphs, two pages suffice up to homeomorphism.

\begin{theorem}
A graph is planar if and only if it is homeomorphic to a graph
of book thickness at most two.
\end{theorem}
\bop Let $G$ be any planar graph and let $G'$ be any subdivision of $G$ such
that all cycles have even length. For instance, one may take the first barycentric subdivision obtained by subdividing each edge exactly once. By Corollary 2, $G'$
is subhamiltonian. \eop

Enomoto and Miyauchi \cite{eno1999} consider “homeomorphic book embeddings”
where edges can cross the spine (i.e., an edge may use more than one page),
and they note that a graph is planar if and only if it has a homeormorphic
book embedding in two pages. It follows from our results above and some
standard topological graph theory arguments that a planar graph with $p$
vertices has a 2-page homeomorphic book embedding with at most $p - 2$
crossings of the spine.


\begin{thebibliography}{99}

\bibitem{gat1968} G. Atneosen, {\bf On the embeddability of compacta in $n$-books: intrinsic and extrinsic properies}, Ph.D. Dissertation, Michigan State University, 1968.

\bibitem{bhky1982}J. Battle, F. Harary, Y. Kodama, and J. W. T. Youngs, Additivity of the genus of a graph, {\it Bull. Amer. Math.} Soc.{\bf 68} (1962),
569-571.

\bibitem{bk79} F. Bernhart and P. C. Kainen, The book thickness of a graph,{\it  J. Combin. Theory Ser. B} {\bf 27} (1979) no. 3, 320–331.


\bibitem{char1971} G. Chartrand, D. P. Geller, and S. Hedetniemi, Graphs with forbidden subgraphs, {\it J. Combin. Theory Ser. B} {\bf 10} (1971) 12-41.

\bibitem{clr1987} F. R. K. Chung, F. T. Leighton, and A. L. Rosenberg, Embeddinggraphs in books: A layout problem with applications to VLSI
design,{\it  SIAM J. on Alg. and Disc. Methods} {\bf 8} (1987) no. 1, 33–58.

\bibitem{dhutch1991} A. M. Dean and J. P. Hutchinson, Relations among embedding parameters for graphs, in {\bf  Graph Theory, Combinatorics,
and Applications}, vol. I, Wiley-Interscience Publ., New York, 1991, pp. 287–296.

\bibitem{dynn1999} I. A. Dynnikov. Three-page link presentation and an untangling algorithm, in {\bf Proc. of the International Conference LowDimensional Topology and Combinatorial Group Theory}, Chelyabinsk, July, 31 - August 7, 1999; Kiev, 2000; pp.112–130.

\bibitem{eno1999} H. Enomoto and M. S. Miyauchi, Embedding graphs into a three page book with {\it O(MlogN)} crossings of edges over the spine,
{\it SIAM J. Discrete Math.} {\bf 12} (1999) 337–341.

\bibitem{gal1989} Z.Galil, R.Kannan and E.Szemeredi. On nontrivial separators for k-page graphs and simulations by nondeterministic one-tape Turing machines, {\it J. Computer and System Sci.}, {\bf 38} (1989) no.1, 134–149.

\bibitem{gle1999} P. M. Gleiss and P. F. Stadler, Relevant cycles in biopolymers and random graphs, 4th Slovene International Conference in Graph
Theory, Bled, June 28 – July 8, 1999.

\bibitem{harary} F. Harary, {\bf Graph Theory}, Addison-Wesley, Reading, MA, 1969.

\bibitem{Has1997} T. Hasunuma and Y. Shibata, Embedding de Bruijn, Kautz and shuffle-exchange networks in books,{\it  Discrete Applied Math.} {\bf 78} (1997) 103–116.

\bibitem{Heath} L. S. Heath and S. Istrail, The pagenumber of genus $g$ graphs is $O(g)$, {\it J. of the ACM} {\bf} 39 (1992) 479–501.


\bibitem{Hol} D. Holton, B. Manvel, B. D. McKay, Hamiltonian cycles in cubic 3-connected bipartite planar graphs, {\it J. Combin. Theory Ser. B} {\bf 38} (1985) 279-297.

\bibitem{PCKa} P. C. Kainen, Some recent results in topological graph theory, in {\bf Graphs and combinatorics}, R. A. Bari and F. Harary, Eds., Springer Lecture Notes in Math. {\bf 406}, Berlin, 1974, pp. 76–108.

\bibitem{PCKb} P. C. Kainen, The book thickness of a graph, II, {\it Congr. Numer.}  {\bf 71} (1990) 127-132.

\bibitem{MWW} D. J. Muder, M. L. Weaver, and D. B. West, Pagenumber of complete bipartite graphs, {\it J. Graph Theory} {\bf 12} (1988) 469–489.

\bibitem{Overbay} S. Overbay, {\bf Generalized book embedings}, Ph. D. Dissertation, Colorado State University, Fort Collins, CO, 1998.

\bibitem{dbw} D. B. West, {\bf Introduction to Graph Theory}, 2nd edition, Prentice-Hall, Upper Saddle River, NJ, 2001.

\bibitem{Whit} H. Whitney, A theorem on graphs,{\it  Ann. of Math.} {\bf 32} (1931) 378–390.

\bibitem{Yann} M. Yannakakis, Linear and book embeddings of graphs, in {\bf VLSI Algorithms and Architectures, Aegean Workshop on
Computing} (Loutraki, Greece, July 8-11, 1986), F. Makedon, K. Mehlhorn, T. S. Papatheodorou, P. G. Spirakis (Eds.), Lecture Notes in Computer Science, vol. 227, Springer 1986.

\end{thebibliography}
\end{document}